\documentclass[a4paper, 10pt, leqno]{article}
\usepackage{pkg}
\usepackage{ncmd}

\usepackage[utf8]{inputenc}

\title{Continuous Cross Approximation of Matrices Arising Out of Kernel Functions\thanks{submitted to the editors \today}}
\author[4,1]{Sumit Singh\,\orcidlink{0009-0002-5581-5349}\thanks{\url{sumit1315singh@gmail.com}, \url{ma22d027@smail.iitm.ac.in}}}
\author[1,2,3,4]{Sivaram Ambikasaran\,\orcidlink{0000-0003-2978-6281}\thanks{\url{sivaambi@dsai.iitm.ac.in}, \url{sivaambi@alumni.stanford.edu}}}

\affil[1]{Wadhwani School of Data Science and Artificial Intelligence, IIT Madras, Chennai, India}
\affil[2]{Robert Bosch Center for Data Science and Artificial Intelligence, IIT Madras, Chennai, India}
\affil[3]{Department of Data Science and Artificial Intelligence, IIT Madras, Chennai, India}
\affil[4]{Department of Mathematics, IIT Madras, Chennai, India}
\date{}

\begin{document}

\maketitle

\begin{abstract}
We propose a residual energy-based framework for constructing low-rank approximations of kernel matrices arising from continuous kernel functions. The method operates in a continuous setting and is based on the adaptive selection of pivot nodes, referred to as \emph{optimal nodes}, which are chosen to minimize the residual energy at each step. This leads to a sequence of rank-$1$ updates of the residual kernel and admits a natural interpretation as a continuous analog of Adaptive Cross Approximation (ACA).
From a theoretical perspective, we show that the residual kernels remain in the class of compact operators and that the approximation error is exactly characterized by the residual energy. We provide convergence guarantees showing that the method yields monotonic error reduction under an alignment condition and achieves geometric decay under practically motivated assumptions. Extensive numerical experiments demonstrate that the proposed method achieves approximation errors close to those of the truncated singular value decomposition across a range of kernel functions. The method exhibits strong robustness with respect to sampling and maintains stable performance across different discretizations. Furthermore, the close agreement between the continuous residual energy and the discrete approximation error highlights the consistency of the formulation. These results establish the proposed approach as a theoretically grounded, practically effective continuous counterpart to classical cross-approximation techniques.
\end{abstract}


\paragraph{MSc Classification:} 65F55, 65W25, 15A23, 47A58, 47G10, 41A45, 65D15.

\section{Introduction}

Kernel functions arise in a wide range of applications in engineering and scientific computing, including partial differential equations \cite{greengard1987fast, massei2022hierarchical}, integral equations \cite{ho2016hierarchical}, inverse problems \cite{ambikasaran2013large}, machine learning \cite{cortes1995support, huang2006extreme}, and Gaussian process models \cite{ambikasaran2015fast}. In many such settings, interactions between domains give rise to kernel matrices. 
More concretely, let $\mclk : \mclx \times \mcly \to \mathbb{R}$ be a kernel function defined over domains $\mclx$ and $\mcly$, and let $X = \{x_1,\dots,x_m\} \subset \mclx$ and $Y = \{y_1,\dots,y_n\} \subset \mcly$ be sets of discretization points. The associated kernel matrix $\boldsymbol{K} \in \mathbb{R}^{m \times n}$ is defined by $\boldsymbol{K}_{ij} = \mclk(x_i, y_j) $, which can be interpreted as $\boldsymbol{K} = \mclk(X, Y)$, and such matrices arise naturally from the discretization of integral equations like
\begin{equation}
    a(x)u(x) + \int_{\mcly} \mclk(x,y)\,u(y)\,dy = f(x) , \quad \forall x \in \mclx
\end{equation}
and the discretization leads to dense linear system $a_iu_i + \sum_j\boldsymbol{K}_{ij}u_j = f_i$.

Although the matrix $\boldsymbol{K}$ is typically dense and usually not low-rank as a whole, the sub-matrix $\mclk(\bar{X}, \bar{Y})$ corresponding to \emph{admissible} subsets $\bar{X}\subset X$ and $\bar{Y}\subset Y$, exhibits low-rank structure \cite{khan2024hodlrdd, singh2025rank}. Hence, efficiently computing a low-rank factorization of such sub-matrices would significantly reduce computational cost.

However, constructing an exact low-rank factorization is generally computationally expensive in practice, as it requires access to the entire kernel matrix and relies on computationally expensive global methods such as the Singular Value Decomposition (SVD), rank-revealing QR factorizations (RRQR), or eigenvalue-based decompositions for symmetric kernels. These approaches typically require $\mclo{mn\min\{m,n\}}$ operations and full matrix access, making them impractical for large-scale problems. Consequently, one seeks an approximate factorization of the form \begin{equation}
    \boldsymbol{K} \approx CUR \quad \text{or} \quad \boldsymbol{K} \approx AB^\top,
\end{equation}
where $C \in \bbr^{m \times r}$, $U \in \bbr^{r \times r}$, $R \in \bbr^{r \times n}$, and $A \in \bbr^{m \times r}, B \in \bbr^{n \times r}$.

\paragraph{Earlier Works:}

A wide range of techniques has been developed to construct such low-rank approximations. Among the available approaches, Adaptive Cross Approximation (ACA), introduced by Bebendorf and Rjasanow \cite{bebendorf2003adaptive}, is an efficient technique widely used for constructing low-rank approximations of kernel matrices. ACA approximates the matrix via successive rank-$1$ updates, selecting pivot rows and columns, and requires access only to the selected rows or columns rather than the full matrix. This makes it particularly attractive for large-scale problems involving smooth kernels.

Despite its practical efficiency, ACA may encounter convergence issues in certain situations \cite{borm2005hybrid}, especially when the pivot selection fails to capture the kernel's dominant components. Moreover, since the method is purely algebraic and operates on a given discretization, it generally lacks rigorous convergence guarantees and may be sensitive to the distribution of sampling points (for instance, see \autoref{fig: boxplot_all_kernels}).

In contrast to purely algebraic approaches, a different class of methods exploits the smoothness of the kernel to construct low-rank approximations via analytic interpolation. \cite{yarvin1998generalized} approximates smooth kernels $\mclk(x,y)$ over $\mclx \times \mcly$ using classical interpolation schemes, which leads to approximations of the form
\begin{equation}
    \mclk(x,y) \approx L_{x}(x,\tilde{X})\,\mclk(\tilde{X},\tilde{Y})^{-1}\,L_{y}(y,\tilde{Y})^\top,
\end{equation}
where $L_x$ and $L_y$ are Lagrange interpolation basis functions defined on $\tilde{X}$ and $\tilde{Y}$ respectively. This can be further recompressed by performing rank-revealing factorizations on $\mclk(\tilde{X},\tilde{Y})$, such as the SVD \cite{fong2009black}. Bebendorf in \cite{bebendorf2000approximation} builds a low-rank factorization of the form \begin{equation} \label{equ: kernel_frctorization_form}
    \mclk(x,y) \approx \mclk(x,\tilde{Y})\,\mclk(\tilde{X},\tilde{Y})^{-1}\,\mclk(\tilde{X},y),
\end{equation}
where  $\tilde{X}$ and $\tilde{Y}$ are interpolation nodes of the interpolant $\mclk(x,y)$ built iteratively. On the other hand, Skeletonized Interpolation (SI) \cite{cambier2019fast} constructs low-rank approximations by first forming the weighted kernel matrix
\begin{equation}
    K_w = \mathrm{diag}(\overline{W}_X)^{1/2} \, \mclk(\overline{X}, \overline{Y}) \, \mathrm{diag}(\overline{W}_Y)^{1/2},
\end{equation}
where $\overline{X}$ and $\overline{Y}$ are tensor-product grids of Chebyshev nodes in $\mclx$ and $\mcly$, respectively, and $\overline{W}_X$ and $\overline{W}_Y$ denote the associated quadrature weights. A subset of nodes $\tilde{X} \subset \overline{X}$ and $\tilde{Y} \subset \overline{Y}$ is then selected via strong rank-revealing QR factorizations applied to $K_w^\top$ and $K_w$, respectively. This procedure yields a low-rank approximation of the form \eqref{equ: kernel_frctorization_form}.

\paragraph{Our Contribution:} 

The methods discussed above highlight a fundamental trade-off in existing approaches for low-rank kernel approximation. Algebraic techniques such as ACA provide adaptivity through pivot selection, but their performance depends on the underlying discretization and may suffer from stability and convergence issues. On the other hand, interpolation-based approaches such as SI exploit the smoothness of the kernel and offer improved stability by selecting nodes from structured grids, but it lacks adaptivity to the specific kernel and domain configuration.

Motivated by these complementary strengths and limitations, we propose, to the best of our knowledge, a new method to determine \emph{``pivots''} by minimizing a global residual energy over continuous domains, thereby providing a \emph{systematic, discretization-independent} criterion for node selection. More specifically, given domains $\mclx$ and $\mcly$, we consider the problem of selecting node sets $S = \{s_1,\dots,s_r\} \subset \mcly$ and $T = \{t_1,\dots,t_r\} \subset \mclx$ that minimize the residual energy functional
\begin{equation}
    E(S,T) = \int_{\mclx} \int_{\mcly} \bigl(\mclk(x,y) - \mclk(x,S)\,\mclk(T,S)^{-1}\,\mclk(T,y)\bigr)^2 \, dx \, dy.
\end{equation}

However, directly minimizing $E(S, T)$ is computationally challenging because the objective depends nonlinearly on the node locations through kernel evaluations and the inverse of $\mclk(T, S)$, yielding a high-dimensional, nonconvex optimization problem. To address this, we adopt a sequential strategy in which the node pairs are selected iteratively. Starting with the initial kernel $K_0 = \mclk$, at each step $k$ we evaluate the current residual kernel $K_{k-1}$ to select a new pair $(t_k, s_k)$ that minimizes the remaining residual energy. We then obtain the updated residual kernel $K_k$ via a rank-$1$ correction. This results in a progressive construction of the node sets $S_k = \{s_1,\dots,s_k\} $ and $ T_k = \{t_1,\dots,t_k\}$, yielding a sequence of low-rank approximations of increasing accuracy.

The resulting nodes are optimal in the sense that they minimize the global residual energy in a greedy manner, and we therefore refer to them as \emph{optimal nodes}. For simplicity, we adopt this terminology throughout the paper, with the understanding that optimality is defined with respect to the chosen energy functional.

The proposed method can be interpreted as a \emph{continuous analog of Adaptive Cross Approximation (ACA)}, which we refer to as \emph{Continuous Cross Approximation (CCA)}. In contrast to classical ACA, which performs pivot selection on a discretized matrix, CCA selects nodes directly in the continuous domain by minimizing a global residual. At each iteration, the method performs a rank-$1$ update of the residual kernel in the Hilbert--Schmidt space $L^2(\mclx \times \mcly)$, subtracting a separable term of the form $K_k(\cdot,s_k)K_k(t_k,\cdot)$. This establishes a direct connection between ACA-type updates and a variational formulation at the operator level.

From this perspective, the proposed method bridges the gap between algebraic approaches, such as ACA, and interpolation-based methods, such as SI. It retains the adaptive nature of ACA while incorporating a continuous, domain-dependent selection criterion for node selection, thereby eliminating dependence on a specific discretization. Moreover, it provides a natural framework for theoretical analysis, as the approximation process can be studied directly in function spaces.

\paragraph{Highlights of the Article:} The main contributions of this work are summarized as follows:
\begin{itemize}
    \item We propose a residual-based framework for low-rank approximation of kernel operators, in which pivots are selected adaptively by minimizing an energy functional in the continuous domain. The resulting pivot nodes, referred to as \emph{optimal nodes}, minimize the residual energy at each iteration and provide an alternative to both algebraic pivoting and fixed interpolation strategies. This leads to a method we term \emph{Continuous Cross Approximation} (CCA), a continuous analog of ACA, in which the approximation is constructed via successive rank-$1$ updates to the residual kernel.

    \item We establish key theoretical properties of the proposed method, including preservation of compactness, an exact characterization of the approximation error, and monotone decay of the residual energy, which corresponds to monotone decay of the error, under the alignment condition. We also provide a spectral interpretation that shows how the operator's dominant components are progressively captured.

    \item Through numerical experiments, we demonstrate that the proposed method achieves approximation errors close to the optimal truncated singular value decomposition while exhibiting strong robustness with respect to sampling. We further show that the loss of optimality at higher ranks is primarily due to numerical instability arising from ill-conditioning, rather than an inherent limitation of the approximation framework. 

    \item Through numerical experiments across multiple kernels and discretizations, we demonstrate that the proposed method achieves approximation errors close to those of the optimal truncated singular value decomposition, is consistently more accurate and stable than ACA across kernels, and also competes well with SI. We further show that deviations from optimality at higher ranks correlate with the growth of the condition number of the interpolation matrix, indicating that the observed degradation in error quality as rank increases is primarily due to numerical instability rather than an inherent limitation of the approximation framework.
     
\end{itemize}

\paragraph{Outline of the Article} The remainder of the paper is organized as follows. \autoref{sec: Preliminaries} introduces the necessary background and notation. \autoref{sec: Energy-Based Node Selection} introduces the algorithm and the implementable version of the algorithm.  \autoref{sec: Theoretical Analysis} presents the theoretical analysis. Numerical experiments are discussed in \autoref{sec: Numerical Experiments}, followed by concluding remarks in \autoref{sec: conclusion}.

\section{Preliminaries} \label{sec: Preliminaries}


Let $(\mclx, dx)$ and $(\mcly, dy)$ be measurable spaces and $L^2(\mclx)$ and $L^2(\mcly)$ be the Hilbert spaces of Lebesgue measurable, square-integrable functions on $\mclx$ and $\mcly$, respectively, equipped with the inner products
\begin{equation} \label{equ: inner product on L^2(X) and L^2(Y)}
    \langle f, g \rangle_{L^2(\Omega)} = \int_{\mclx} f(x) g(x)\,dx,\quad \text{where } \Omega = \{ \mclx, \mcly \} 
\end{equation}
and corresponding norms $\|f\|_{L^2(\mclx)} = \langle f,f \rangle_{L^2(\mclx)}^{1/2}$ 
and $\|g\|_{L^2(\mcly)} = \langle g,g \rangle_{L^2(\mcly)}^{1/2}$. Both spaces are complete, hence Hilbert spaces; see \cite{hackbusch2012integral}. We further define the product space $\mclh := L^2(\mclx \times \mcly),$ equipped with the inner product
\begin{equation}
    \langle K_1, K_2 \rangle_{\mclh} = \int_{\mclx} \int_{\mcly} K_1(x,y)\, K_2(x,y)\, dx\, dy,
\end{equation}
and norm $\|K\|_{\mclh} = \langle K,K \rangle_{\mclh}^{1/2}$.

Now given a kernel $\mclk \in \mclh$, it defines a bounded integral operator \cite[Thm.~3.2.7]{hackbusch2012integral} $ \mclt : L^2(\mcly)\to L^2(\mclx) $ as 
\begin{equation}\label{equ: integral operator}
    (\mclt \phi)(x) = \int_{\boldsymbol{\mathcal{Y}}} \mclk(x,y) \phi(y)\,dy.
\end{equation} 
The operator $\mclt$ depends on the kernel $\mclk$, where the operator $\mclt$ is called a \emph{Hilbert-Schmidt integral operator} and the kernel $\mclk$ is called \emph{Hilbert-Schmidt kernel}. This framework includes most kernels used in practice, such as Green's functions, smooth kernels, radial basis function (RBF) kernels, and logarithmic kernels away from singularities. 

Now we take a look at a well-known theory of the kernel function $\mclk$, which is one of the article's backbones. Given a continuous bivariate function $\mclk:\mclh\to \mclr$, we can define the Singular Value Expansion (SVE) \cite[\S~5]{schmidt1907theorie} as \begin{equation} \label{equ: singular value expression of mclk}
    \mclk(x,y) = \sum_{k=1}^\infty \sigma_k \phi_k(x)\psi_k(y)
\end{equation}
where $\sigma_k$'s are \emph{singular values} and $\phi_k\in L^2(\mclx)$, $\psi_k\in L^2(\mcly)$ are referred to as singular functions which are orthogonal with respect to the suitable inner product \eqref{equ: inner product on L^2(X) and L^2(Y)} in $L^2(\mclx)$ and $L^2(\mcly)$ respectively \cite[\S~5]{schmidt1907theorie}. The singular values form a real, non-increasing sequence. It can be shown that the only limit point of the singular values for square integrable functions is 0 \cite[\S~5]{schmidt1907theorie}\footnote{Schmidt refers to singular values as eigenvalues as his work primarily dealt with symmetric kernels, and the modern notion of singular values had not been developed when this work was published.}. Each term in \eqref{equ: singular value expression of mclk} is an ``\emph{outer product}" of two univariate functions $\phi_k $ and $ \psi_k$, called a \emph{rank-1 function}. 

The sequence of singular values decays rapidly for smooth functions. The smoother the underlying function, the more rapid the decay of its singular values. This behavior is discussed in several literature; for instance, see \cite{little1984eigenvalues, griebel2010approximation, hackbusch2012integral, wendland1995piecewise}.


The SVE \eqref{equ: singular value expression of mclk} naturally leads to finite-rank approximations of the kernel. For any $r \in \mathbb{N}$, the truncated expansion
\begin{equation} \label{equ: truncated sve}
    \mclk_r(x,y) := \sum_{k=1}^r \sigma_k \phi_k(x)\psi_k(y)
\end{equation}
defines a rank-$r$ approximation of $\mclk$. Now, a classical result states that $\mclk_r$ is the best rank-$r$ approximation of $\mclk$ in the $\mclh$-norm, i.e., $\|\mclk - \mclk_r\|_{\mclh} = \inf_{\mathrm{rank}(K_r)\le r} \|\mclk - K_r\|_{\mclh}$ \cite[\S~2]{townsend2013extension}\cite{hackbusch2012integral}. Moreover, the approximation error is given by
\begin{equation} \label{equ: sve error decay}
    \|\mclk - \mclk_r\|_{\mclh}
    = \bkt{\sum_{k=r+1}^\infty \sigma_k^2}^{\frac{1}{2}},
\end{equation}
which characterizes the decay of the optimal approximation error. Even though \eqref{equ: truncated sve} provides an optimal approximation, its computation requires access to the singular functions and values of the associated integral operator, which is typically infeasible in large-scale or data-driven settings. In practice, one therefore considers structured rank-$r$ approximations of the form
\begin{equation} \label{equ: separable form}
    K_r(x,y) = \sum_{k=1}^r u_k(x)\, v_k(y),
\end{equation}
where $u_k \in L^2(\mclx)$ and $v_k \in L^2(\mcly)$. Such representations arise naturally in numerical methods for kernel approximation \cite[\S\S~1.3]{martinsson2025fast}, and the goal of these methods is to construct structured approximations $K_r$ whose error is close to the optimal error $\|\mclk - \mclk_r\|_{\mclh}$.

A particularly important class of such approximations is given by interpolation-based or skeleton-type representations. Given the sets of pivots nodes $T \subset \mclx$, $S \subset \mcly$, we construct approximations of the form
\begin{equation} \label{equ: skeleton form}
    K_r(x,y) = \mclk(x,S)\, \mclk(T,S)^{-1}\, \mclk(T,y)
\end{equation}
Here, $\mclk(x, S)$ and $\mclk(T, y)$ denote vectors of kernel evaluations. Such representations arise in various low-rank approximation techniques, including cross approximation and skeletonization methods \cite{cambier2019fast, bebendorf2000approximation, hackbusch2012integral}.

The formulation in \eqref{equ: skeleton form} reduces the kernel approximation problem to selecting suitable node sets $T$ and $S$. The design of efficient strategies for selecting these nodes is a central challenge in low-rank kernel approximation, and in \autoref{sec: Energy-Based Node Selection} we propose a method for generating such a set of nodes. Before moving further, let us define the terminology that we use.


\paragraph{Well-separated domains:}
    Let $\mclx, \mcly \subset \bbr^d$ be bounded domains. We say that $\mclx$ and $\mcly$ are \emph{well-separated} if there exists a constant $\eta > 0$ such that $\operatorname{dist}(\mclx,\mcly) \;\ge\; \eta \,\max\{\operatorname{diam}(\mclx),\,\operatorname{diam}(\mcly)\}$, where $\operatorname{dist}(\mclx,\mcly) = \inf_{x\in\mclx,\;y\in\mcly} \|x-y\|, $ and $ \operatorname{diam}(\mclx) = \sup_{x,x'\in\mclx} \|x-x'\|$.

It is well known that kernel functions restricted to well-separated domains exhibit low-rank structure. This phenomenon is closely related to the smoothness and analytic behavior of the kernel away from the diagonal. The following classical result from approximation theory shows that the analytically extendable functions can be efficiently approximated by separable low-rank expansions. 

\begin{mylma}[Polynomial approximation of analytic functions] \label{lma: polynomial approximation of analytic functions}
Let $f: V = [-1,1]^d \subset \mathbb{R}^d \to \mathbb{R}$ be analytic and admit an analytic extension to a generalized Bernstein ellipse $\mathcal{B}(V,\rho)$, where $\rho > 1$. Suppose that 
\begin{equation}
    \|f\|_{\infty^*} = \max_{y \in \mathcal{B}(V,\rho)} |f(y)| \le M.
\end{equation}
Then, its multivariate polynomial interpolant $\tilde{f}$ satisfies
\begin{equation}
    \|f - \tilde{f}\|_{\infty^*} \le \frac{4 M V_d}{\rho - 1}\,\rho^{-p},
\end{equation}
where $p$ is a predefined constant and $V_d$ is a constant depending on the dimension $d$ and $\rho$.
\end{mylma}

The \autoref{lma: polynomial approximation of analytic functions} is discussed in \cite{khan2024hodlrdd}. Also, a detailed discussion on the generalized Bernstein ellipse and analytic continuation can be found there. The generalized version of \autoref{lma: polynomial approximation of analytic functions} is stated and proved in \cite{glau2019improved}. 

The above \autoref{lma: polynomial approximation of analytic functions} provides an error bound in the supremum norm over a Bernstein-type domain and implies exponential convergence for analytic functions. In the context of kernel approximation, such estimates justify approximating $\mclk(x,y)$ by interpolating in one variable while treating the other as a parameter, leading to separable representations of the form
\begin{equation}
    \mclk(x,y) \approx \sum_{k=1}^r \mclk(x,y_k)\,L_k(y),
\end{equation}
where $\{y_k\}$ are interpolation nodes and $L_k$ are the corresponding basis functions.

Although the error estimates are stated in the supremum norm, they imply corresponding decay in the $L^2$ norm on bounded domains \cite[\S 3]{rudin1987real}. This observation justifies the use of low-rank representations of the form \eqref{equ: separable form} and factorizations based on interpolation nodes. While classical approaches determine these nodes either through fixed interpolation grids or algebraic pivoting strategies, our objective is to select them adaptively by minimizing a global residual energy in the continuous setting. This perspective serves as the basis for the framework developed in this work.

\section{Energy-Based Node Selection} \label{sec: Energy-Based Node Selection}


As discussed in the previous section, low-rank approximations of kernels depend on structured representations, such as separable expansions or interpolation-based forms. The effectiveness of such approximations depends on the selection of suitable node sets. In this section, we develop an energy-based framework for constructing these nodes by iteratively minimizing a global residual energy, leading to a sequence of progressively improved low-rank approximations.

Let $K_0 = \mclk$ and recursively define, for each $k \ge 1$, the residual energy functional \begin{equation}\label{equ: residual functional at kth srep}
    E_k(t,s) = \int_\mclx\int_\mcly \bkt{K_{k-1}(x,y) - \frac{K_{k-1}(x,s)K_{k-1}(t,y)}{K_{k-1}(t,s)} }^2 dxdy
\end{equation}
This functional quantifies the global $\mclh$-error obtained after subtracting a rank-$1$ cross approximation of the residual kernel $K_{k-1}$ based on the pair $(t,s)$. Now, at each iteration, we select the nodes $ (t_k,s_k) = \arg\min\limits_{t\in\mclx,s\in\mcly} E_k(t,s) $ and accordingly define the updated residual kernel by 
\begin{equation}\label{equ: kernel updating}
    K_{k}(x,y) = K_{k-1}(x,y) - \frac{K_{k-1}(x,s_k)K_{k-1}(t_k,y)}{K_{k-1}(t_k,s_k)}, \quad \text{ with } K_{k-1}(t_k,s_k) \neq 0.
\end{equation}

This recursive construction generates a sequence of \emph{residual kernels} $\{K_k\}_{k\ge 0}$ with decreasing energy, reflecting successive removal of rank-$1$ components. This procedure can be viewed as a continuous analog of adaptive cross approximation, in which pivot selection is performed via global energy minimization as described in \autoref{alg: the ideal algorithm}.
\begin{algorithm}[ht]
    \begin{algorithmic}[1]
        \State \textbf{Input:} Initial kernel function $\mathcal{K}(x,y)$, desired rank $r$, domains $ \boldsymbol{\mathcal{X}} $ and $ \boldsymbol{\mathcal{Y}} $
        \State \textbf{Output:} Optimal pair of nodes  $\{ (t_k, s_k) \}_{k=1}^r$
        \State Initialize $K_0 \gets \mathcal{K}$ \Comment{Start with the original kernel}
        \For{$k = 1$ to $r$}
            \State Find the pivot pair $(t_k, s_k)$: \[ (t_k,s_k) = \arg\min_{\substack{t\in \boldsymbol{\mathcal{X}} \\ s \in \boldsymbol{\mathcal{Y}}}}E_k(t,s) 
            \]
            \State Update kernel:
            \[ K_k(x,y) \gets K_{k-1}(x,y) - \frac{K_{k-1}(x,s_k)\, K_{k-1}(t_k,y)}{K_{k-1}(t_k,s_k)} \]
        \EndFor
    \end{algorithmic}
    \caption{Optimal Node Selection via Residual Updates in the ideal situation}
    \label{alg: the ideal algorithm}
\end{algorithm}

While \autoref{alg: the ideal algorithm} provides a conceptually clear formulation of optimal pivot selection through global energy minimization, in practice, its direct implementation is computationally infeasible. The reason behind this is that the kernel $\mclk$ itself and the residual kernels $K_k(x,y)$ might not admit a unique global minimizer. In general, the associated objective $E_k$ may contain multiple local minima due to the structure of the kernel and the geometry of the domains $\mclx$ and $\mcly$, making it a highly nonconvex optimization problem over the continuous domains. 

To overcome this limitation, we adopt a practical strategy. Instead of solving the global minimization problem exactly, we generate multiple initial guesses and perform local optimization to identify candidate pivot pairs. These candidates are then evaluated using a surrogate of the residual energy, computed via matrix approximations on sampled sets $X \subset \mclx$ and $Y \subset \mcly$.

This leads to the \autoref{alg: optimal-rank-r-residual-update1}, an implementable version of the \autoref{alg: the ideal algorithm}, which can be interpreted as a discrete counterpart of the ideal continuous algorithm.

\begin{algorithm}[ht]
\caption{Optimal Node Selection for Rank-$r$ Approximation of Kernel via Residual Updates}
\label{alg: optimal-rank-r-residual-update1}
\begin{algorithmic}[1]
\State \textbf{Input:} Initial kernel function $\mathcal{K}(x,y)$, desired rank $r$, domains $x \in \boldsymbol{\mathcal{X}} $, $y \in \boldsymbol{\mathcal{Y}} $
\State \textbf{Output:} Optimal pair of nodes  $\{ (t_k, s_k) \}_{k=1}^r$

\State Initialize $K \gets \mathcal{K}$ \Comment{Start with the original kernel}

\State Construct the kernel matrix $ \boldsymbol{K} = \mclk(X,Y)  $ \Comment{ $X \subset \boldsymbol{\mathcal{X}}, Y \subset \boldsymbol{\mathcal{Y}}$}

\For{$k = 1$ to $r$}
    \State Generate multiple initial guesses for $(t, s)$ in the domain
    \For{each guess}
        \State Find a candidate pair $(t, s)$ \[ (t,s) = \arg\min_{\substack{t\in \boldsymbol{\mathcal{X}} \\ s \in \boldsymbol{\mathcal{Y}}}}\int\limits_{\boldsymbol{\mathcal{X}}}\int\limits_{\boldsymbol{\mathcal{Y}}}\bkt{K(x,y) - \frac{K(x,s)\, K(t,y)}{K(t,s)}}^2dxdy \]
        \State Construct a temporary pivot set $S=\{s_1,s_2,\dotsc,s_{k-1},s\}$, $T=\{t_1,t_2,\dotsc,t_{k-1},t\}$ to get the approximated kernel matrix        
        \[ \widehat{\boldsymbol{K}}^{(k)} = \mclk(X,S)\mclk(T,S)^{-1}\mclk(T,Y)  \] 
        \State {Evaluate the candidate by computing the relative approximation error:}
        \[
        \boldsymbol{\varepsilon}^{(k)}(t,s) = \frac{ \| \boldsymbol{K} - \widehat{\boldsymbol{K}}^{(k)} \|_F }{ \| \boldsymbol{K} \|_F }
        \]
    \EndFor
    \State Find \[ (t_k, s_k) = \arg\min_{t,s}\boldsymbol{\varepsilon}^{(k)}(t,s) \] 
    \State Update kernel:
    \[
        K(x,y) \gets K(x,y) - \frac{K(x,s_k)\, K(t_k,y)}{K(t_k,s_k)}
    \]
\EndFor
\end{algorithmic}
\end{algorithm}

Although \autoref{alg: the ideal algorithm} is motivated by a global energy minimization principle, its performance is not evident from the construction alone. In the next section, we provide a rigorous analysis of the method. In particular, we establish structural properties of the residual kernels, the energy-decay property, and structural representations of the kernel matrices based on them.

\section{Theoretical Analysis} \label{sec: Theoretical Analysis}

In this section, we analyze the proposed residual-based kernel approximation procedure. Our goal is to understand the structural and analytical properties of the sequence of residual kernels $\{K_k\}_{k\geq 0}$ generated by \autoref{alg: the ideal algorithm}, as well as the behavior of the associated approximation error. Now, motivated by practical implementations where the domains are chosen as hypercubes in $\bbr^d$, we henceforth assume that $\mclx, \mcly \subset \bbr^d$ are compact sets.

We begin by interpreting the construction from an operator-theoretic perspective, which allows us to leverage classical results from compact operator theory and provides a natural framework for the analysis. In particular, each residual kernel $K_k$ induces a linear integral operator $\mclt_k : L^2(\mcly) \to L^2(\mclx)$ defined by
\begin{equation}
    (\mclt_k \phi)(x) = \int_\mcly K_k(x,y)\,\phi(y)\,dy.
\end{equation}
We first show that the recursive construction preserves fundamental regularity and compactness properties, ensuring that the approximation process remains within the class of compact operators, a key feature for establishing convergence and low-rank structure.

\begin{mylma}
Let the initial kernel $K_0(x,y)$ be continuous on $\mclx \times \mcly$, where $\mclx$ and $\mcly$ are compact sets. Then, for every $k \ge 0$, the residual kernels $K_k(x,y)$ remain continuous on $\mclx \times \mcly$, and the associated operators $\mclt_k$ are compact.
\end{mylma}

\begin{proof}
We prove the continuity of $K_k$ by induction. Since $K_0(x,y)$ is continuous on $\mclx\times\mcly$, the statement holds for $k=0$. Now, assume that $K_{k-1}(x,y)$ is continuous. Then from the update formula given in \autoref{equ: residual functional at kth srep} we observe that $x\mapsto K_{k-1}(x,s_k)$ and $y\mapsto K_{k-1}(t_k,y)$ are continuous functions. Therefore the product $K_{k-1}(x,s_k) K_{k-1}(t_k,y)$ is continuous on $\mclx\times\mcly$. Since subtraction preserves continuity, $K_k(x,y)$ is continuous. Thus, by induction, $K_k\in C(\mclx\times\mcly)$ for all $k$. Since $\mclx$ and $\mcly$ are compact, an integral operator with continuous kernel defines a compact operator from $L^2(\mcly)$ to $L^2(\mclx)$. Therefore, each $\mclt_k$ is compact.
\end{proof}

Now, as each residual kernel $K_k$ defines a compact integral operator $\mclt_k$, we move next to quantify the approximation process. 
To measure the residual, we therefore introduce the \emph{residual energy} associated with $K_k$ and study its behavior under the recursive update. In particular, we show that the energy admits an exact characterization in terms of the residual functional $E_k(t,s)$ and exhibits a monotone decay under a suitable alignment condition.

\begin{mylma}\label{lma: energy decay}
    Let the residual energy be defined using the updated kernel as defined by \autoref{equ: kernel updating} as \begin{equation}
        \mcle_k =  \int_\mclx\int_\mcly K_k(x,y)^2 dx dy,\qquad \text{ for } k = 1,2,\dotsc
    \end{equation} 
    then the residual energy satisfies \begin{enumerate}[(i)]
        \item $\mcle_{k}=E_k(t_k, s_k)$
        \item $\mcle_{k+1} \leq \mcle_{k}$, if f the selected pair $(t_k,s_k)$ satisfies $2K_k(t_k,s_k) \langle \mclt_kv_{k}, u_{k} \rangle \geq  \magn{u_{k}}^2\magn{v_{k}}^2 $.
    \end{enumerate}
\end{mylma}

\begin{proof}
    By definition of residual energy, we have \begin{equation}
        \mcle_{k} =  \int_\mclx\int_\mcly K_{k}(x,y)^2 dx dy.
    \end{equation}
    \begin{enumerate}[(i)]
        \item Substituting $K_{k}$ using \autoref{equ: kernel updating} in the above expression, we have \begin{equation}
            \mcle_{k} =  \int_\mclx\int_\mcly \bkt{K_{k-1}(x,y) - \frac{K_{k-1}(x,s_k)K_{k-1}(t_k,y)}{K_{k-1}(t_k,s_k)} }^2 dxdy
        \end{equation} 
        But the right-hand side is exactly the definition of $E_k(s_k,t_k)$. Therefore \[ \mcle_{k}=E_k(t_k, s_k) \]

        \item By construction of the \autoref{alg: optimal-rank-r-residual-update1}, we have $ (t_k,s_k) = \arg\min\limits_{t\in\mclx,s\in\mcly} E_k(t,s) $. Therefore, \begin{equation}
            E_k(t_k,s_k)\leq E_k(t,s)\quad \forall \; (t,s)\in \mclx\times\mcly
        \end{equation} 
        Now, from part (i) we have $\mcle_{k+1}\leq E_{k+1}(t,s)$ $\forall \; (t,s)$ and hence it will also hold for $(t_k,s_k)$. Now, using \autoref{equ: residual functional at kth srep}, we have \begin{align}
            \mcle_{k+1}\leq E_{k+1}(t_k,s_k) 
            &= \mcle_k - \frac{2 \langle \mclt_kv_{k}, u_{k} \rangle_{L^2(\mclx)} }{K_k(k,k)} + \frac{ \magn{u_{k}}^2\magn{v_{k}}^2 }{K_k(k,k)^2}
        \end{align}
        where, $ u_{k}(x) = K_k(x,s_k) $ and $ v_{k}(y) = K_k(t_k,y)$. Thus we have $\mcle_{k+1} \leq \mcle_{k}$, provided the following holds \begin{equation}\label{equ: condition on monotone decay}
            2K_k(t_k,s_k) \langle \mclt_kv_{k}, u_{k} \rangle \geq  \magn{u_{k}}^2\magn{v_{k}}^2 .
        \end{equation}
    \end{enumerate}      
\end{proof}

\begin{myremark}
    The above result shows that the algorithm performs an optimal rank-$1$ correction at each step with respect to the residual energy functional. This is analogous to \emph{weak greedy approximation schemes}, where each iteration selects the best element from a restricted dictionary to reduce the residual.
\end{myremark}

The above formulation captures the monotone decay of the residual in energy form. However, it remains unclear how each iteration modifies the dominant singular components of $\mclt_k$. To address this, we now proceed to a spectral analysis of the residual operators. By viewing $\mclt_k$ as a Hilbert-Schmidt operator, we can express the residual energy in terms of its singular values. This provides a way to track the successive removal of the dominant singular components.

Since $K_k\in \mclh$, the operator $\mclt_k$ is a Hilbert-Schmidt operator, and its Hilbert-Schmidt norm satisfies \begin{equation}
    \magn{\mclt_k}_{HS}^2 = \int_\mclx\int_\mcly K_k(x,y)^2\,dxdy = \mcle_k.
\end{equation} Thus, the residual energy $\mcle_k$ corresponds exactly to the squared Hilbert-Schmidt norm of the operator $\mclt_k$.

The kernel update formula given in \autoref{equ: residual functional at kth srep} corresponds to subtracting a rank-$1$ operator from $\mclt_{k-1}$, assuming $\mclt_0 = \mclt $. Now, defining $u_k(x)=K_{k-1}(x,s_k),$ and $v_k(y)=K_{k-1}(t_k,y)/K_{k-1}(t_k,s_k) ,$ we can write $\mclt_k = \mclt_{k-1} - u_k\otimes v_k$, where the following denotes a rank-$1$ operator \begin{equation}
    (u_k\otimes v_k)\phi = u_k \int_\mcly v_k(y)\phi(y)\,dy.
\end{equation}

Now, since the operator $\mclt_{k-1}$ is compact, by the Hilbert-Schmidt spectral theorem, it admits a singular value decomposition which is given in  kernel form below
\begin{equation}
    K_{k-1}(x,y) = \sum_{j=1}^{\infty} \sigma_j^{(k-1)} u_j^{(k-1)}(x) v_j^{(k-1)}(y),
\end{equation} 
where $\{\sigma_j^{(k-1)}\}$ are nonnegative singular values and $\{u_j^{(k-1)}\}$, $\{v_j^{(k-1)}\}$ are orthonormal systems in $L^2(\mclx)$ and $L^2(\mcly)$, respectively. In this representation, the residual energy satisfies 
\begin{equation}
    \mcle_{k-1} = \sum_{j=1}^{\infty} \big(\sigma_j^{(k-1)}\big)^2.
\end{equation}

Now, the rank-$1$ structure of the update allows us to control the evolution of the singular values across iterations. Since $\mclt_k = \mclt_{k-1} - u_k \otimes v_k$ can be interpreted as a rank-$1$ perturbation of $\mclt_{k-1}$, Weyl's interlacing inequality \cite[Chapter II]{gohberg1978introduction} for singular values of compact operators implies that, for all $j \geq 1$,
\begin{equation}\label{equ: interlacing inequality}
    \sigma_{j+1}^{(k-1)} \;\le\; \sigma_j^{(k)} \;\le\; \sigma_{j-1}^{(k-1)}.
\end{equation}

Thus, each rank-$1$ update to the residual can locally shift the singular values at most one position. In particular, no new large singular values can be introduced, and the dominant part of the spectrum can only decrease in a controlled manner. As a consequence, the residual energy $\mcle_k$ 
evolves through a sequence of singular values that remain interlaced with those of the previous iteration. This ensures a form of \emph{controlled spectral decay} of the residual energy, i.e., each iteration removes at most one significant singular component, while the remaining spectrum is only mildly perturbed. This provides a spectral explanation for the gradual reduction of the residual energy.


\begin{myremark}
    
While Weyl's interlacing inequality \eqref{equ: interlacing inequality} ensures that the singular value cannot increase under a rank-$1$ update, it does not by itself guarantee strict decay. In particular, one only has $\sigma_2^{(k-1)} \le \sigma_1^{(k)} \le \sigma_1^{(k-1)},$ and therefore a reduction in the leading singular value occurs only when the rank-$1$ correction is sufficiently aligned with the dominant singular components.
\end{myremark}

Therefore, subtraction of rank-$1$ terms in the update step may be interpreted as removing the dominant low-rank component of the residual operator if the alignment condition \eqref{equ: condition on monotone decay} holds, which is discussed in the following \autoref{rmk: singular direction by alignment condition}.

\begin{myremark}\label{rmk: singular direction by alignment condition}
The leading singular vectors maximize $\langle \mathcal{T}_k v, u \rangle $ over all $\magn{u}=\magn{v}=1$ and therefore they represent the directions of strongest interaction. The selection $(u_k, v_k)$ does not enforce this optimality. Instead, the condition \eqref{equ: condition on monotone decay} ensures $\langle \mathcal{T}_k v_k, u_k \rangle$ remains sufficiently large. Consequently, while $(u_k, v_k)$ need not coincide with the leading singular vectors, 
they must have nontrivial components along the leading singular directions of $\mathcal{T}_k$, and hence, it captures a significant portion of the operator’s dominant singular directions, explaining the monotone decay behavior of the residual energy.
\end{myremark}

This interpretation suggests that the proposed procedure approximates the operator $\mclt_0$ by successively removing rank-$1$ components, analogous to truncated singular-value expansions. Although the monotone decay condition \eqref{equ: condition on monotone decay} is not imposed explicitly in \autoref{alg: the ideal algorithm}, the minimization of the residual functional $E_k(t,s)$ naturally favors pairs $(t,s)$ for which the inner product $\langle \mclt_k v_t, u_s \rangle$ is relatively large compared to $\|u_s\|^2\|v_t\|^2/K_k(t,s)$.

\begin{figure}[H]
    \centering
    
    \begin{subfigure}[t]{0.48\linewidth}
        \centering
        \includegraphics[width=\linewidth]{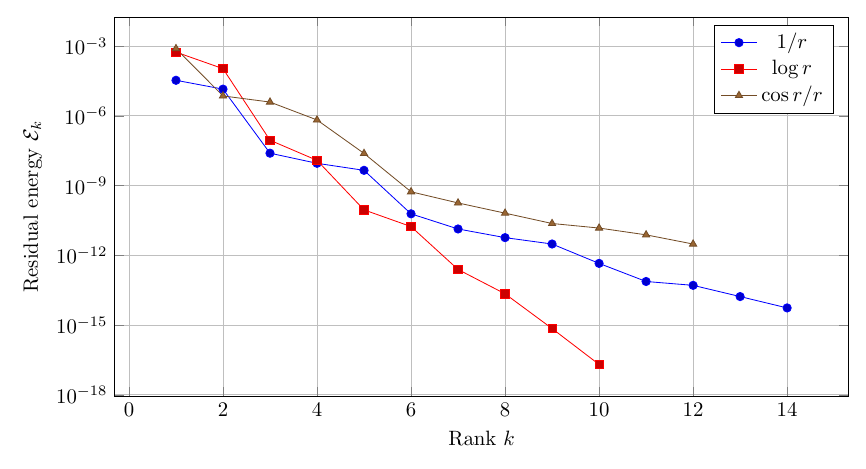}
        \caption{$\mclx=[0,1] \times [0,1] $ and $\mcly=[2,3] \times [2,3] $.}
    \end{subfigure}
    \hfill
    \begin{subfigure}[t]{0.48\linewidth}
        \centering
        \includegraphics[width=\linewidth]{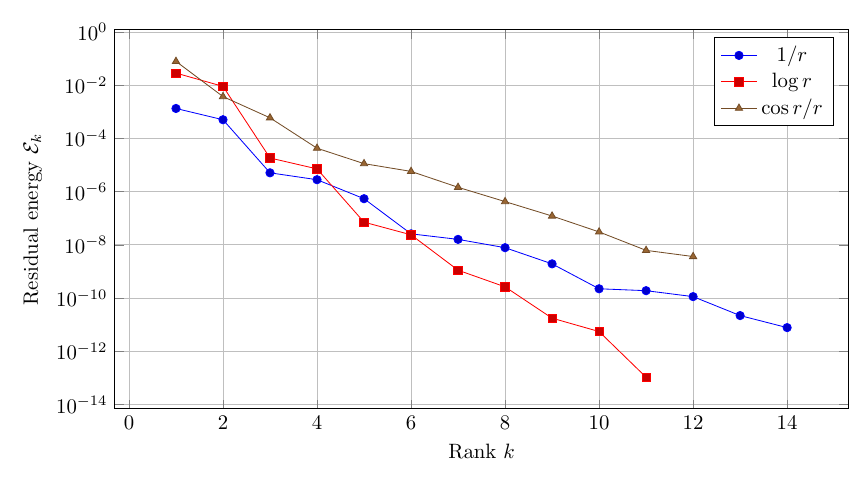}
        \caption{$\mclx=[-3,-1] \times [0,2] $ and $\mcly=[1,3] \times [0,2] $.}
    \end{subfigure}
    
    \caption{ Residual energy decay for the kernels $1/r$, $\log r$, and $\cos r/r$ over different well-separated domains. The plots illustrate the rapid decrease in residual energy as the rank increases.}
    
    \label{fig: residual energy decay in 2D}
\end{figure}

Consequently, even without explicitly enforcing \eqref{equ: condition on monotone decay}, the optimization procedure tends to select rank-one corrections that remove a significant portion of the residual energy. This explains why monotone decay of the residual energy is consistently observed in numerical experiments (see \autoref{fig: residual energy decay in 2D}).

Now, if each update removes a fixed proportion of the residual energy, then under an assumption observed in practice, one can expect a geometric decay of the sequence of residual kernels. This is formalized in the following result.

\begin{myremark}[Geometric convergence of residual kernel]
\label{rmk: geometric hs convergence}
Let $(K_k)_{k\ge0}\subset \mclh$ be the sequence of residual kernels generated by \autoref{alg: the ideal algorithm}, and suppose there exists a constant $c\in(0,1)$ such that $ \mcle_k-\mcle_{k+1} \ge c\,\mcle_k $ for every $k\ge0$, then $K_k \to 0$ in $\mclh$-norm at a geometric rate. 
\end{myremark}

\begin{proof} Rearranging the assumption, we have $\mcle_{k+1}\le (1-c)\mcle_k$. Applying this recursively, we obtain $\mcle_k \leq (1-c)^k\mcle_0 $. Since $0<c<1$, we have $(1-c)^k\to0$ as $k\to \infty$, and hence $\mcle_k\to0$. Since $\mcle_k=\magn{\mclt_k}_{HS}^2$ it follows that $\magn{\mclt_k}_{HS}\to0$. Now, since $K_k \in \mclh$ and $\mclh$ is isometrically isomorphic\footnote{The mapping $K_k \mapsto \mclt_k$, defined by $(\mclt_k \phi)(x) = \int_{\mcly} K_k(x,y)\phi(y)\,dy$, is an isometric isomorphism between $\mclh$ and the space of Hilbert–Schmidt operators from $L^2(\mcly)$ to $L^2(\mclx)$.} to the space of Hilbert–Schmidt operators, we have $\|\mclt_k\|_{HS} = \|K_k\|_{\mclh}.$ Therefore, $\|K_k\|_{\mclh} \to 0$, which proves the result.
\end{proof}

\begin{figure}[ht]
    \centering
    
    \begin{subfigure}[t]{0.48\linewidth}
        \centering
        \includegraphics[width=\linewidth]{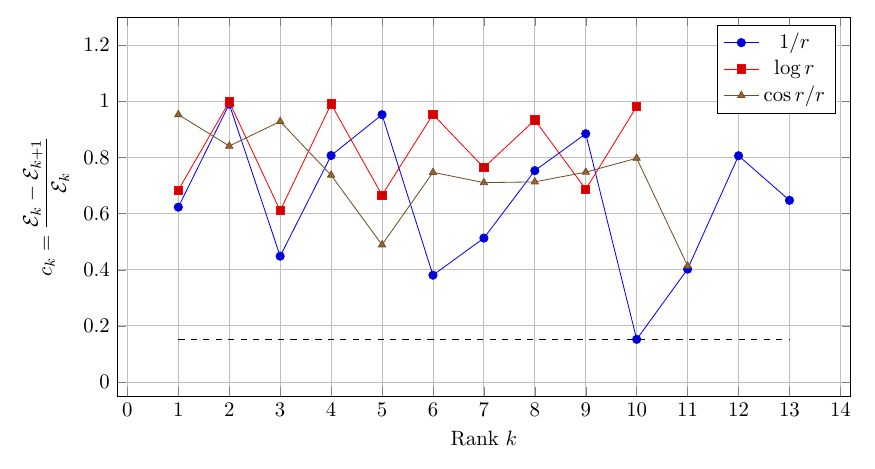}
        \caption{ $\mclx=[-3,-1]\times[0,2]$ and $\mcly =[1,3]\times[0,2]$.}
 
    \end{subfigure}
    \hfill
    \begin{subfigure}[t]{0.48\linewidth}
        \centering
        \includegraphics[width=\linewidth]{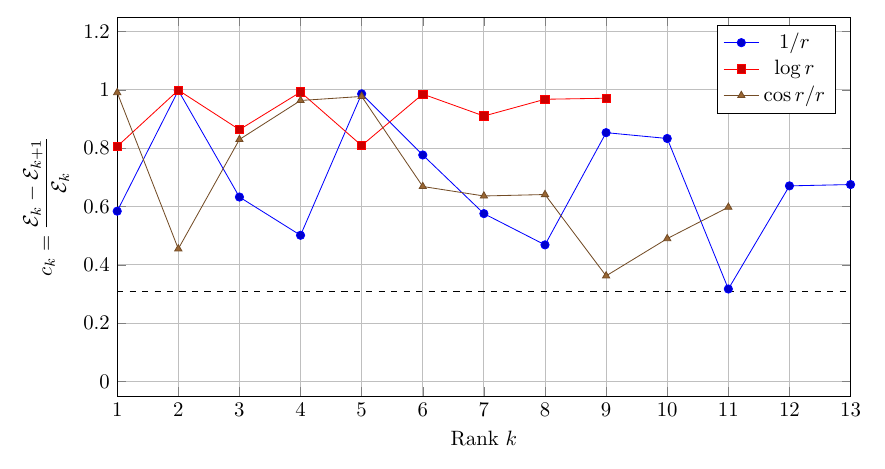}
        \caption{$\mclx=[0,1]\times[0,1]$ and $\mcly=[2,3]\times[2,3]$.}
    \end{subfigure}
    
    \caption{The observed values of $c_k$ for different kernels remain uniformly bounded away from zero in the two-dimensional setting. This indicates consistent decay of the residual energy across successive ranks.}
    
    \label{fig: c_k values in 2D}
\end{figure}

The assumption in \autoref{rmk: geometric hs convergence} is motivated by consistent observations across our numerical experiments. We analyze the empirical decay ratio $c_k = \frac{\mcle_k - \mcle_{k+1}}{\mcle_k}.$ \autoref{fig: c_k values in 2D} shows the behavior of $(c_k)_k$ for different kernel functions and domain configurations in the two-dimensional setting. We observe consistently across our numerical experiments that although the sequence $(c_k)_k$ exhibits mild fluctuation, it remains consistently bounded away from zero, i.e., there exists a constant $c_* > 0$ such that $c_k \geq c_*$ for all $k$  and shows no downward trend as $k$ increases.

These observations support our assumption that each iteration removes a fixed proportion from the residual energy, and thereby justifying the geometric convergence used in the analysis.

The geometric decay in the residual energy suggests that the iterative updates effectively capture the kernel's dominant components. Now, the following result shows that the approximation admits an explicit separable representation in terms of the selected pivot nodes. 

\begin{mylma} \label{lma: iterative_cross_form}
Let $\mclk\in \mclh $ be the kernel function, and let $\{(t_i,s_i)\}_{i=1}^k \subset \mclx \times \mcly$ be the pivot nodes selected by \autoref{alg: the ideal algorithm}. Then \begin{equation}
    \widehat{\mclk}_k(x,y) = \sum_{j=0}^{k-1} \frac{K_j(x,s_{j+1})K_j(t_{j+1},y)} {K_j(t_{j+1},s_{j+1})}
\end{equation}
be the rank-$k$ approximation of $\mclk$ obtained from the recursive residual updates as defined in \eqref{equ: kernel updating}. Then, for every $k\ge 1$, we have \begin{enumerate}[(i)]
    \item $\mclk(x,y) = \widehat{\mclk}_k(x,y) + K_k(x,y)$,
    \item $\widehat{\mclk}_k(x,y)=\mclk(x,S_k)\,\mclk(T_k,S_k)^{-1}\,\mclk(T_k,y),\; \text{ where } S_k=\{s_1,\dots,s_k\}, T_k=\{t_1,\dots,t_k\}.$
\end{enumerate} 
\end{mylma}

\begin{proof} We proceed by induction on $k$ for both cases.
\begin{enumerate}[(i)]
    \item For $k=1$, we have from \eqref{equ: kernel updating}, $K_1(x,y) = \mclk(x,y) - \frac{\mclk(x,s_1)\mclk(t_1,y)}{\mclk(t_1,s_1)}$. Rearranging, we obtain\begin{equation}
        \mclk(x,y) = \frac{\mclk(x,s_1)\mclk(t_1,y)}{\mclk(t_1,s_1)} + K_1(x,y) = \widehat{\mclk}_1(x,y)+K_1(x,y).
    \end{equation}
    
    Now let us assume that $\mclk(x,y)=\widehat{\mclk}_k(x,y)+K_k(x,y)$ is true for some $k\ge 1$. Now, rearranging the residual kernel update formula \eqref{equ: kernel updating}, we have 
    \begin{equation}
        K_k(x,y) = \frac{K_k(x,s_{k+1})K_k(t_{k+1},y)} {K_k(t_{k+1},s_{k+1})} + K_{k+1}(x,y).
    \end{equation}
    
    Substituting this into the induction hypothesis yields $\mclk(x,y)=\widehat{\mclk}_{k+1}(x,y)+K_{k+1}(x,y)$.
    
    \item For $k=1$, we have $\widehat{\mclk}_1(x,y) = {\mclk(x,s_1)\mclk(t_1,y)}/{\mclk(t_1,s_1)}$. Now, as $\mclk(x,S_1)=\mclk(x,s_1), \mclk(T_1,S_1)=\mclk(t_1,s_1), \mclk(T_1,y)=\mclk(t_1,y)$, it follows immediately that $\widehat{\mclk}_1(x,y)=\mclk(x,S_1)\,\mclk(T_1,S_1)^{-1}\,\mclk(T_1,y).$

    Assume now that the result holds for some $k\ge 1$, i.e., $\widehat{\mclk}_k(x,y)=\mclk(x,S_k)\,M_k^{-1}\,\mclk(T_k,y),$ where $M_k:=\mclk(T_k,S_k)$.
    We now show that the identity holds for $k+1$.
    
    Let us define $M_{k+1} = \mclk(T_{k+1},S_{k+1}) = \begin{bmatrix} M_k & a\\ b & c \end{bmatrix}$,  where $a=\mclk(T_k,s_{k+1}), b=\mclk(t_{k+1},S_k), c=\mclk(t_{k+1},s_{k+1}) $.  Using the block inverse formula, \begin{equation}
        M_{k+1}^{-1} = \begin{bmatrix} M_k^{-1}+M_k^{-1}a\Delta^{-1}bM_k^{-1} & -M_k^{-1}a\Delta^{-1} \\[1mm] -\Delta^{-1}bM_k^{-1} & \Delta^{-1} \end{bmatrix},\text{ where } \Delta=c-bM_k^{-1}a.
    \end{equation}
    Now, by the induction hypothesis, $\widehat{\mclk}_k(t_{k+1},s_{k+1}) = bM_k^{-1}a$, and hence $\Delta = \mclk(t_{k+1},s_{k+1})-\widehat{\mclk}_k(t_{k+1},s_{k+1}) = K_k(t_{k+1},s_{k+1}).$
    Now writing  $\mclk(x,S_{k+1}) = \begin{bmatrix} \mclk(x, S_k) & \mclk(x,s_{k+1}) \end{bmatrix}$, and $\mclk(T_{k+1},y) = \begin{bmatrix} \mclk(T_k,y)\\ \mclk(t_{k+1},y) \end{bmatrix}$ and substituting the block inverse expression and simplifying, we obtain \begin{align}
    \mclk(x,S_{k+1})M_{k+1}^{-1}\mclk(T_{k+1},y) = \widehat{\mclk}_k(x,y) + \frac{ \bkt{\mclk(x,s_{k+1})-\mclk(x,S_k)M_k^{-1}a} \bkt{\mclk(t_{k+1},y)-bM_k^{-1}\mclk(T_k,y)} }{\Delta}.
    \end{align}
    Using (i), we have $K_k(x,s_{k+1}) = \mclk(x,s_{k+1})-\widehat{\mclk}_k(x,s_{k+1}) $, and $K_k(t_{k+1},y) = \mclk(t_{k+1},y)-\widehat{\mclk}_k(t_{k+1},y). $ Hence, using the induction hypothesis, we have  \begin{equation}
        \mclk(x,S_{k+1})M_{k+1}^{-1}\mclk(T_{k+1},y) = \widehat{\mclk}_k(x,y) + \frac{ K_k(x,s_{k+1})K_k(t_{k+1},y) }{ K_k(t_{k+1},s_{k+1}) }.
    \end{equation}
    Hence, we have \begin{equation}
        \widehat{\mclk}_{k+1}(x,y) = \mclk(x,S_{k+1})\mclk(T_{k+1},S_{k+1})^{-1}\mclk(T_{k+1},y).
    \end{equation}
\end{enumerate}
This completes the induction.
\end{proof}

In summary, the analysis in this section establishes that the \autoref{alg: the ideal algorithm} generates a sequence of compact operators with monotonically decreasing energy, admits a precise spectral interpretation through rank-$1$ updates, and yields an explicit separable representation in terms of the selected pivot, the \emph{optimal nodes}. In particular, the method can be viewed as a \emph{continuous analog of cross approximation}, where the kernel is progressively reconstructed using the obtained optimal nodes from the domains. 

The theoretical results provide insight into spectral behavior and convergence under suitable conditions of the proposed \autoref{alg: the ideal algorithm}, we now move to the numerical experiments in the section to complement the theoretical findings and to assess the effectiveness of the method in realistic settings. 

\section{Numerical Experiments} \label{sec: Numerical Experiments}

In this section, we present numerical experiments to assess the effectiveness of the proposed method. The theoretical foundation is provided by the ideal residual-update procedure described in \autoref{alg: the ideal algorithm}, where pivot pairs are selected by minimizing the residual energy functional over the continuous domains. This procedure is implemented using the computational framework described in \autoref{alg: optimal-rank-r-residual-update1}, which approximates the continuous optimization via numerical quadrature and multi-start optimization. The experiments are implemented in Julia \cite{bezanson2017julia}. The computational framework follows a sequential greedy residual-update procedure, where each pivot pair $(t_k, s_k)$ is selected by minimizing the residual functional $E_k(t,s)$ associated with the updated kernel.

The evaluation of $E_k(t,s)$ involves approximating a double integral over $\mclx \times \mcly$, for which we use Gauss– Legendre quadrature in each dimension. To solve the resulting nonconvex optimization problem, we use the \texttt{Optim.jl} package with a multi-start strategy, where multiple initial guesses are used to improve robustness. 
The source code for all experiments is available at \url{https://github.com/SAFRAN-LAB/CCA-optimal_nodes/}.

\begin{figure}[ht]
    \centering
    \includegraphics[width=0.4\linewidth]{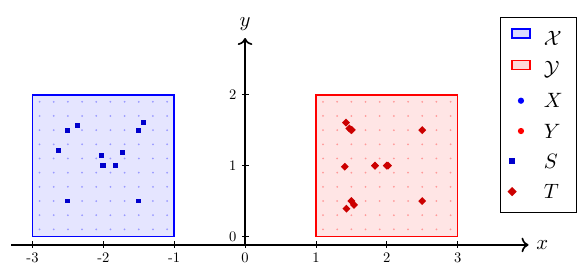}
    \caption{Geometric configuration of the domains. The domains $\mclx = [-3,-1]\times[0,2]$ (blue) and $\mcly = [1,3]\times[0,2]$ (red) are discretized into point sets $X$ and $Y$ respectively. The sets $S$ and $T$ denote the optimal nodes selected by the \autoref{alg: optimal-rank-r-residual-update1} for the kernel $1/r$.}
    \label{fig: domain details}
\end{figure}

We begin with an elementary discussion of the domains, as depicted on \autoref{fig: domain details}. We consider the usual kernels \begin{equation}\label{equ: kernel list}
    \mclk(x,y) = \kappa(\magn{x-y}), \quad \text{ and }\quad \kappa(r)\in \bkct{ \frac{1}{r}, \log r, \frac{\cos r}{r}, \frac{1}{\sqrt{1+r}}, \sqrt{1+r}},
\end{equation}
where  $x \in \mclx , y\in \mcly$. $\mclx$ and $\mcly$ are two squares of side length 2, centered at $(-2, 1)$ and $(2, 1)$ respectively, i.e., $\mclx = [-3, -1]\times[0, 2] $ and $\mcly = [1, 3]\times[0, 2] $. For error evaluation, we construct discretizations of size $65 \times 65$ on each domain, resulting in $n=4225$ sample points per domain. Both uniform grids and randomly sampled points are considered.  Let $X$ and $Y$ denote the discretization points in $\mclx$ and $\mcly$, and let $T$ and $S$ denote the selected pivot nodes. Then, the approximated matrix of the kernel matrix $\boldsymbol{K} =\mclk(X,Y) $ is obtained as  \begin{equation} \label{equ: approximated kernel matrix}
    \widehat{\boldsymbol{K}} = \mclk(X,S)\mclk(T,S)^{-1}\mclk(T,Y). 
\end{equation} 
The relative error due to the approximation is measured in Frobenius norm as ${\|\boldsymbol{K}-\boldsymbol{ \widehat{ \boldsymbol{K}} }\|_F}/{\magn{\boldsymbol{K}}_F}$, and for the random case, we report the mean relative error over $500$ independent trials.

We compare the proposed method with other standard low-rank approximation techniques for kernel matrices, including partially pivoted Adaptive Cross Approximation (ppACA) \cite{bebendorf2003adaptive}, Skeletonized Interpolation (SI) \cite{cambier2019fast}, and the truncated Singular Value Decomposition (SVD). All the methods are evaluated in the same domain settings, and the results are summarized in \autoref{tab: error data table with comparison}. The comparison is done in the following way:
\begin{enumerate}
    \item For a given kernel in a certain domain $X$, and $Y$, obtain the original kernel matrix $\boldsymbol{K} = \mclk(X,Y)$ and for a prescribed target rank $r$ of $\boldsymbol{K}$, obtain the optimal nodes using \autoref{alg: optimal-rank-r-residual-update1}.
    \item Now construct the approximated matrix using the optimal nodes in the same domain as given in \autoref{equ: approximated kernel matrix} and compute the relative error with $\boldsymbol{K}$.
    \item For the same rank and kernel in the domain $X$ and $Y$, build the ppACA approximation and evaluate its relative error in Frobenius norm with $\boldsymbol{K}$.
    \item For Skeletonized Interpolation, select \emph{pivot nodes} from Chebyshev grids on $\mclx$ and $\mcly$, construct the corresponding approximation, and then compute the relative error.
\end{enumerate}

\begin{table}[ht]
\centering
\begin{tabular}{llcccc}
\toprule
 &  & \multicolumn{2}{c}{Uniform Grid} & \multicolumn{2}{c}{Random Grid} \\
\cmidrule(lr){3-4} \cmidrule(lr){5-6}
Kernel & Method & Rank & Relative Error & Rank & Mean Relative Error \\
\midrule

\multirow{3}{*}{$1/r$}
 & Optimal Nodes & 14 & $3.160371\times10^{-6}$ & 14 & $2.603404\times10^{-6}$ \\
 & SVD           & 14 & $2.510962\times10^{-7}$ & 14 & $1.829823\times10^{-7}$ \\
 & ppACA         & 14 & $9.953475\times10^{-6}$ & 14 & $8.775568\times10^{-6}$ \\
 & SI            & 14 & $1.975709\times10^{-6}$ & 14 & $1.569133\times10^{-6}$ \\
\midrule

\multirow{3}{*}{$\log r$}
 & Optimal Nodes & 11 & $8.451296\times10^{-8}$ & 11 & $6.635491\times10^{-8}$ \\
 & SVD           & 11 & $1.647485\times10^{-8}$ & 11 & $1.176472\times10^{-8}$ \\
 & ppACA         & 11 & $3.080435\times10^{-7}$ & 11 & $3.425584\times10^{-7}$ \\
 & SI            & 11 & $8.447236\times10^{-8}$ & 11 & $7.070352\times10^{-8}$  \\
\midrule

\multirow{3}{*}{$\frac{\cos r}{r}$}
 & Optimal Nodes & 12 & $7.238291\times10^{-5}$ & 12 & $6.117635\times10^{-5}$ \\
 & SVD           & 12 & $1.212579\times10^{-5}$ & 12 & $9.184520\times10^{-6}$ \\
 & ppACA         & 12 & $4.742501\times10^{-4}$ & 12 & $2.035438\times10^{-4}$ \\
 & SI            & 12 & $2.551607\times10^{-5}$ & 12 & $2.255897\times10^{-5}$  \\
 \midrule

\multirow{3}{*}{$\frac{1}{\sqrt{1+r}}$}
 & Optimal Nodes & 9 & $5.949229\times10^{-6}$ & 9 & $5.297518\times10^{-6}$ \\
 & SVD           & 9 & $9.681600\times10^{-7}$ & 9 & $7.588193\times10^{-7}$ \\
 & ppACA         & 9 & $1.583490\times10^{-5}$ & 9 & $1.578853\times10^{-5}$ \\
 & SI            & 9 & $3.579787\times10^{-6}$ & 9 & $3.291803\times10^{-6}$  \\
  \midrule

\multirow{3}{*}{${\sqrt{1+r}}$}
 & Optimal Nodes & 9 & $6.390571\times10^{-6}$ & 9 & $5.797219\times10^{-6}$ \\
 & SVD           & 9 & $6.977305\times10^{-7}$ & 9 & $5.738864\times10^{-7}$ \\
 & ppACA         & 9 & $9.849246\times10^{-6}$ & 9 & $7.139570\times10^{-6}$ \\
 & SI            & 9 & $4.974543\times10^{-6}$ & 9 & $4.397744\times10^{-6}$  \\
\bottomrule
\end{tabular}
\caption{Comparison of relative errors for uniform and random grids across different kernels and approximation methods for the domain $\mclx = [-3, -1]\times[0, 2] $ and $\mcly = [1, 3]\times[0, 2] $.}
\label{tab: error data table with comparison}
\end{table}

\begin{figure}[ht]
    \centering
    \includegraphics[width=0.325\linewidth]{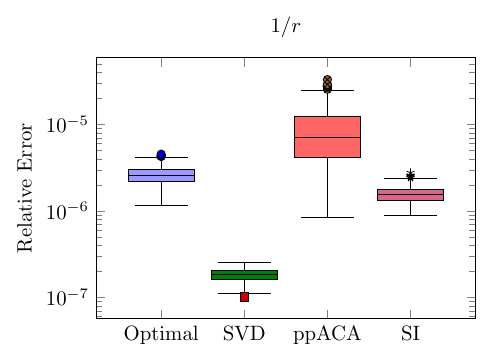} \hfill
    \includegraphics[width=0.325\linewidth]{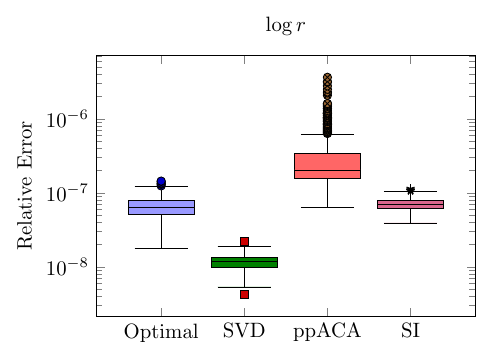} \hfill
    \includegraphics[width=0.325\linewidth]{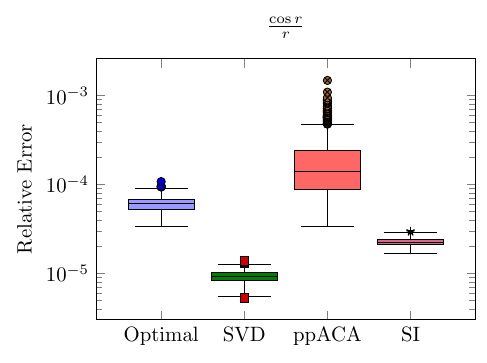} \hfill
    \includegraphics[width=0.325\linewidth]{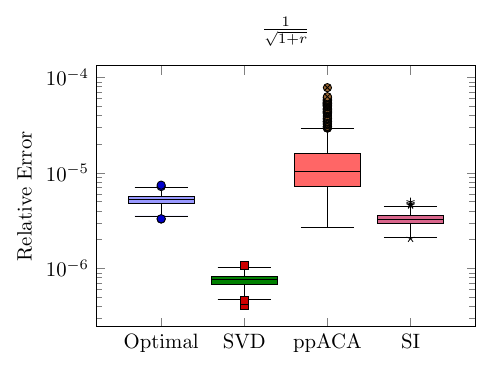}
    \includegraphics[width=0.325\linewidth]{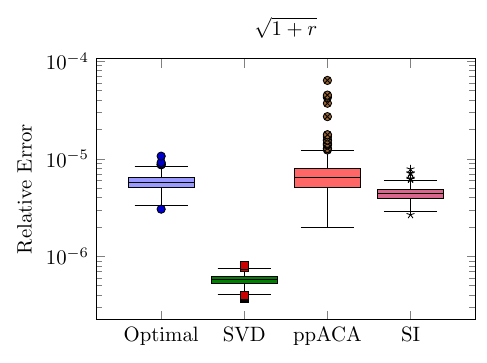}
    \caption{Relative errors computed over $500$ independent random discretizations of the domains $\mclx = [-3,-1]\times[0,2]$ and $\mcly = [1,3]\times[0,2]$ for all the kernels as mentioned in \eqref{equ: kernel list}. For each method and kernel, the distribution of errors across random grids is summarized using boxplots: the central line indicates the median, the box represents the interquartile range (25th to 75th percentiles), whiskers indicate the range of typical (non-outlier) values, and individual points denote outliers.} 
    \label{fig: boxplot_all_kernels}
\end{figure}

To assess the robustness of our method with respect to sampling, we complement the mean-error results reported in \autoref{tab: error data table with comparison} by analyzing the variability of the approximation error across multiple random discretizations of the domains. Even though the \autoref{tab: error data table with comparison} reports the average performance over $500$ trials, it does not capture how sensitive the method is to changes in the sampling. To address this, for each kernel, we generate $500$ independent random grids on $\mclx$ and $\mcly$, compute the corresponding kernel approximations, and record the relative errors. The resulting distributions are summarized using boxplots in \autoref{fig: boxplot_all_kernels}.

A closer look at \autoref{fig: boxplot_all_kernels} reveals consistent patterns across all kernels. The SVD results are tightly concentrated, indicating nearly identical errors across different discretizations and the best possible approximation. The relative error due to the optimal nodes also exhibits a relatively compact inter-quartile range, with only mild spread and few outliers, suggesting that the selected pivots adapt consistently to changes in the sampling. In contrast, the ppACA method shows a notably wider spread, with several extreme values, suggesting that its performance can vary significantly across grid realizations. The Skeletonized Interpolation (SI) method shows a narrow spread, similar to SVD, and our method competes well with it. 

These observations indicate that the proposed method maintains stable performance across different kernels and discretizations. Further improvements in the optimization strategy, such as enhanced initialization or more robust solvers within \autoref{alg: optimal-rank-r-residual-update1}, may lead to additional gains in accuracy and potentially improved performance.

We now turn to a quantitative comparison of our method's approximation error as a function of rank. In particular, we investigate how our method compares with the optimal rank-$k$ approximation given by the truncated singular value decomposition. The residual formulation developed in the \autoref{sec: Theoretical Analysis} provides a natural way to quantify the approximation error at each iteration. In particular, by part (i) of \autoref{lma: iterative_cross_form}, after $k$ steps, the residual kernel $K_k = \mclk - \widehat{\mclk}_k$ represents the exact approximation error and therefore, the squared approximation error in $\mclh$-norm is given by $\mcle_k = \|\mclk-\widehat{\mclk}_k\|_{\mclh}^2 = \|K_k\|_{\mclh}^2$. Now, the squared rank-$k$ truncation error obtained from the singular value decomposition is $\mcle_k^{\mathrm{svd}} = (\sum_{j>k}\sigma_j^2) $.  Since the truncated SVD is the best among all rank-$k$ approximations, it follows that $\mcle_k^{\mathrm{svd}} \leq \mcle_k $.

Now for numerical evaluation, we consider the kernel matrix $ \boldsymbol{K} \in \mathbb{R}^{N\times N}$ over the domain $\mclx \times \mcly$, and compute its singular value decomposition $\boldsymbol{K} = U \Sigma V^\top$. The best rank-$k$ approximation is given by $\boldsymbol{K}_k^{\mathrm{svd}} = U_{:,1:k}\Sigma_{1:k,1:k}V_{:,1:k}^\top$ and the corresponding approximation error is given by $E_k^{\mathrm{svd}} = \| \boldsymbol{K} -\boldsymbol{K}_k^{\mathrm{svd}}\|_F/\magn{\boldsymbol{K}}_F$. Now, using the optimal nodes obtained from \autoref{alg: optimal-rank-r-residual-update1} for all the kernels, we evaluate the error as $ E_k^{\mathrm{opt}} = \| \boldsymbol{K} - \widehat{\boldsymbol{K}}_k \|_F/\magn{\boldsymbol{K}}_F $, where $ \widehat{\boldsymbol{K}}_k $ is obtained by \eqref{equ: approximated kernel matrix}. We also compare the quantities $\Bar{\mcle}_k =  \sqrt{\mcle_k/\mcle_0}$, where $\mcle_0 = \magn{\mclk}_\mclh^2$. 


\begin{figure}[ht]
    \centering  
    \begin{subfigure}[t]{0.32\linewidth}
        \centering
        \includegraphics[width=\linewidth]{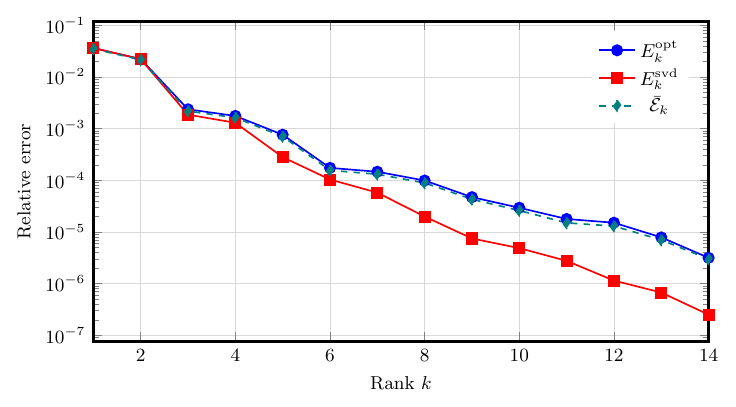}
        \caption{$1/r$}
        \label{fig:error_k1_r14}
    \end{subfigure}
    \hfill
    \begin{subfigure}[t]{0.32\linewidth}
        \centering
        \includegraphics[width=\linewidth]{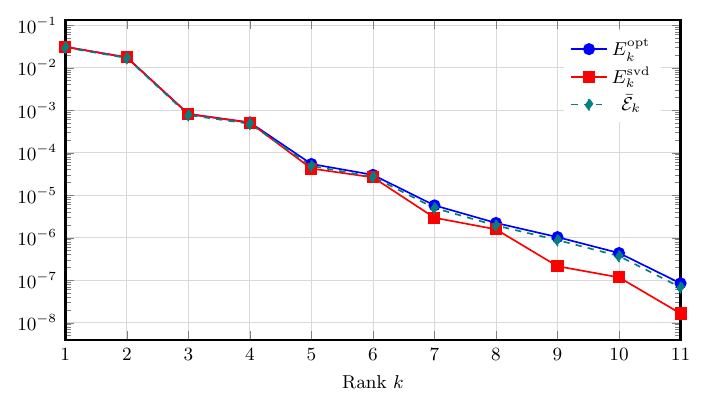}
        \caption{$\log r$}
        \label{fig:error_k2_r11}
    \end{subfigure}
    \hfill
    \begin{subfigure}[t]{0.32\linewidth}
        \centering
        \includegraphics[width=\linewidth]{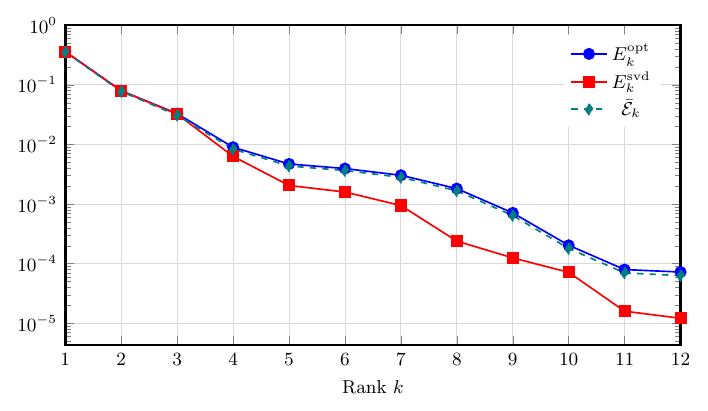}
        \caption{$\frac{\cos r}{r}$}
        \label{fig:error_k4_r12}
    \end{subfigure}
    \hfill
    \begin{subfigure}[t]{0.32\linewidth}
        \centering
        \includegraphics[width=\linewidth]{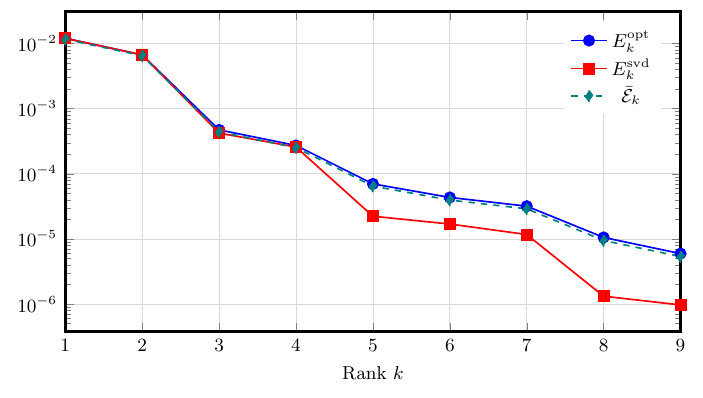}
        \caption{$\frac{1}{\sqrt{1+r}}$}
        \label{fig:error_k5_r09}
    \end{subfigure}
    \begin{subfigure}[t]{0.32\linewidth}
        \centering
        \includegraphics[width=\linewidth]{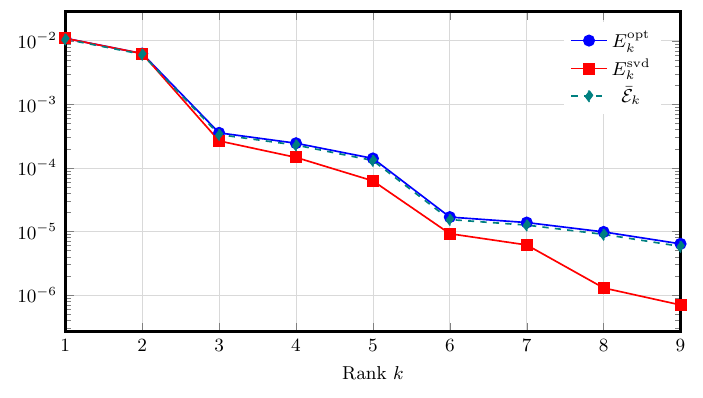}
        \caption{${\sqrt{1+r}}$}
        \label{fig:error_k7_r09}
    \end{subfigure}
    \caption{Comparison of approximation errors as a function of rank $k$ for different kernels. For each kernel, we plot the relative error of our method $E_k^{\mathrm{opt}}$, the optimal truncated SVD error $E_k^{\mathrm{svd}}$, and the normalized residual energy $\bar{\mcle}_k = \sqrt{\mcle_k/\mcle_0}$.}
    \label{fig: error_comparison}
\end{figure}

The error curves in \autoref{fig: error_comparison} exhibit a consistent decay with increasing rank for all kernels. As expected, the SVD-based approximation achieves the lowest error at each rank. The error obtained from our method closely follows the decay, with only a small gap across all kernels. Furthermore, the normalized residual energy $\bar{\mcle}_k$ tracks the discrete approximation error $E_k^{\mathrm{opt}}$ remarkably well, indicating that the continuous residual formulation provides an accurate estimate of the actual approximation error. 

While the error curves indicate that the proposed method closely follows the optimal decay, it is natural to further quantify the gap between the two approximations. To quantify this gap, we consider the ratio $\rho_k = {E_k^{\mathrm{opt}}}/{E_k^{\mathrm{svd}}}$ together with the condition number $\kappa_k = \mathrm{cond}\bkt{\mclk(T_k, S_k)}$. A value of $\rho_k$ close to $1$ indicates near-optimal approximation, while a uniformly bounded ratio suggests quasi-optimal behavior of the approximation.

\begin{figure}[H]
    \centering
    \begin{subfigure}[t]{0.32\linewidth}
        \centering
        \includegraphics[width=\linewidth]{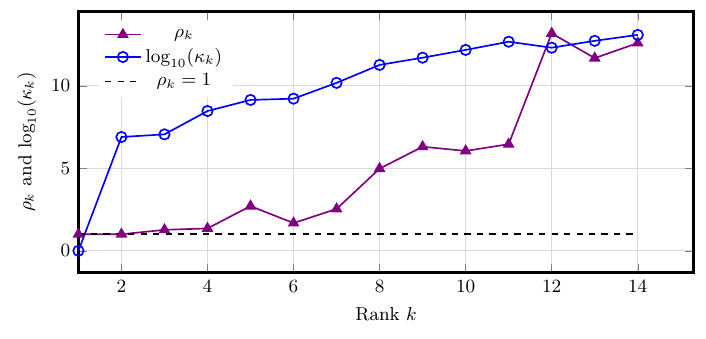}
        \caption{$1/r$}
        \label{fig: rho_k1_r14}
    \end{subfigure}
    \hfill
    \begin{subfigure}[t]{0.32\linewidth}
        \centering
        \includegraphics[width=\linewidth]{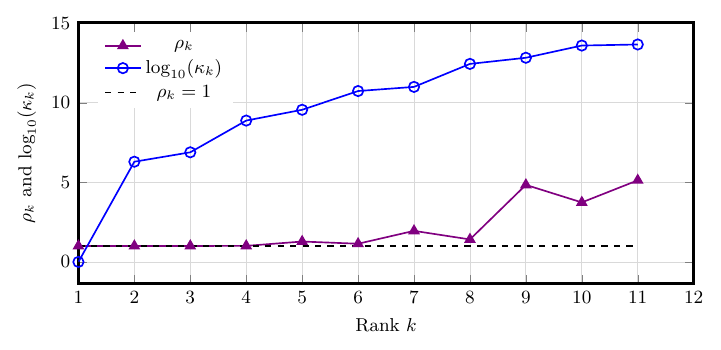}
        \caption{$\log r$}
        \label{fig: rho_k2_r11}
    \end{subfigure}
    \hfill
    \begin{subfigure}[t]{0.32\linewidth}
        \centering
        \includegraphics[width=\linewidth]{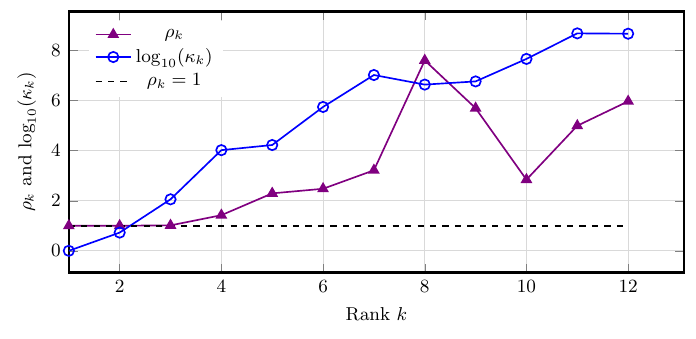}
        \caption{$\frac{\cos r}{r}$}
        \label{fig: rho_k4_r12}
    \end{subfigure}
    \hfill
    \begin{subfigure}[t]{0.32\linewidth}
        \centering
        \includegraphics[width=\linewidth]{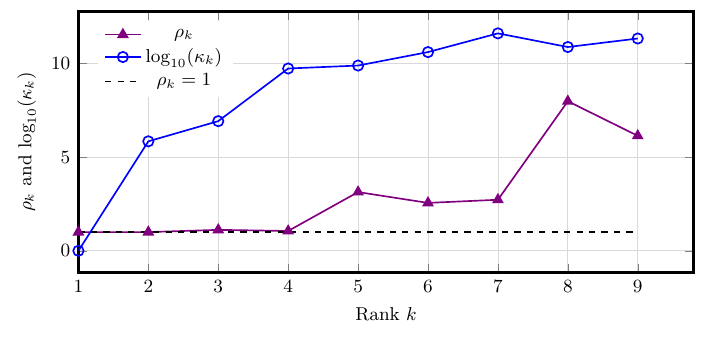}
        \caption{$\frac{1}{\sqrt{1+r}}$}
        \label{fig: rho_k5_r09}
    \end{subfigure}
    \begin{subfigure}[t]{0.32\linewidth}
        \centering
        \includegraphics[width=\linewidth]{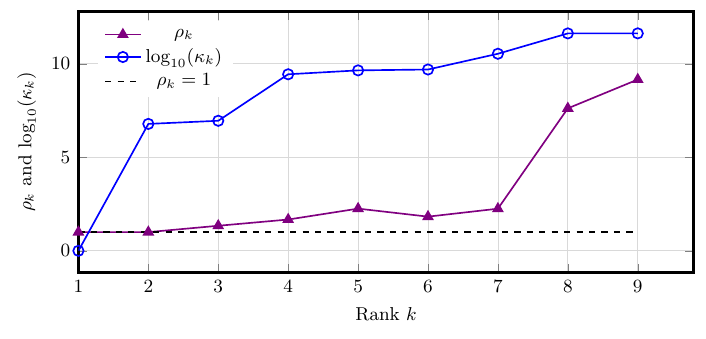}
        \caption{${\sqrt{1+r}}$}
        \label{fig: rho_k7_r09}
    \end{subfigure}
    \caption{Comparison of the approximation ratio $\rho_k = E_k^{\mathrm{opt}} / E_k^{\mathrm{svd}}$ (purple) and the logarithm of the condition number $\log_{10}(\kappa_k)$ (blue) as functions of the rank $k$ for different kernels. The dashed line corresponds to $\rho_k = 1$, which represents the optimal approximation. }
    \label{fig: rho_condition}
\end{figure}

The results in \autoref{fig: rho_condition} show that $\rho_k \approx 1$ at low ranks across all kernels, confirming that our method achieves near-optimal accuracy in this regime. As the rank increases, however, $\rho_k$ exhibits a gradual growth, indicating a deviation from optimality. This trend is closely aligned with the rapid increase in the condition number $\kappa_k$ of the interpolation matrix $\mclk(T_k, S_k)$. In particular, the growth in $\rho_k$ becomes pronounced in the regime where $\kappa_k$ becomes large, often reaching values between $10^8$ and $10^{12}$.

These observations suggest that the proposed method remains quasi-optimal in exact arithmetic, and that the observed degradation at higher ranks is primarily due to numerical instability arising from ill-conditioning of $\mclk(T_k, S_k)$, rather than an inherent limitation of the approximation strategy.

\section{Conclusion}\label{sec: conclusion}

In this work, we proposed a residual-based framework for constructing low-rank approximations of kernel operators by adaptively selecting \emph{pivot nodes}, which we refer to as \emph{optimal nodes}, since each selection minimizes the residual energy. The method is formulated in a continuous setting and can be naturally interpreted as a \emph{continuous analog of Adaptive Cross Approximation (ACA)}, where the approximation is built via successive rank-$1$ updates to the residual kernel. This perspective provides a principled mechanism for node selection and establishes a direct link between the approximation error and the residual energy. We refer to the method \emph{Continuous Cross Approximation (CCA)}.

From a theoretical standpoint, we showed that the sequence of residual kernels remains within the class of compact operators and that the approximation error admits an exact characterization through the residual. Under a uniform energy reduction assumption, we established geometric convergence in the Hilbert–Schmidt norm. In addition, a spectral interpretation based on rank-$1$ perturbations clarifies how the dominant components of the operator are progressively extracted through the iterative updates.

The numerical experiments demonstrate that the proposed method achieves approximation errors close to those of the truncated singular value decomposition across various kernels. In particular, the method exhibits strong robustness with respect to sampling and maintains consistent performance across different discretizations. The close alignment between the continuous residual energy and the discrete approximation error provides strong empirical evidence that the residual energy-based formulation accurately captures the behavior of the underlying matrix approximation problem.

Overall, our proposed method offers a continuous, theoretically grounded alternative to classical cross-approximation methods, providing a unified perspective that bridges operator-level analysis and matrix-level computation. Future directions include the development of stabilization techniques to control ill-conditioning at higher ranks, and integration with hierarchical and fast algorithms for large-scale applications.

\appendix

\section{Equivalence of the recursive and cross-interpolatory forms for \texorpdfstring{$k=2$}{k=2} }
\label{app:k2_equivalence}

In this section, we verify the equivalence between the rank-$2$ approximation and its compact cross-interpolatory representation as given in \autoref{lma: iterative_cross_form}. 

For $k=2$, the recursive approximation is given by
\begin{equation}
    \widehat{\mclk}_2(x,y) = \frac{\mclk(x,s_1)\mclk(t_1,y)}{\mclk(t_1,s_1)} + \frac{K_1(x,s_2)K_1(t_2,y)}{K_1(t_2,s_2)},\; \text{where } K_1(x,y) = \mclk(x,y) - \frac{\mclk(x,s_1)\mclk(t_1,y)}{\mclk(t_1,s_1)}.
\end{equation}
Let $a=\mclk(t_1,s_1), b=\mclk(t_1,s_2), c=\mclk(t_2,s_1), d=\mclk(t_2,s_2)$. Then $ \mclk(T,S)= \begin{bmatrix} a & b\\ c & d \end{bmatrix}. $  Now, we have $K_1(x,s_2)=\mclk(x,s_2)-\frac{b}{a}\mclk(x,s_1), K_1(t_2,y)=\mclk(t_2,y)-\frac{c}{a}\mclk(t_1,y), $ and $ K_1(t_2,s_2)=d-\frac{bc}{a}=\frac{ad-bc}{a} $. Now, substituting these expressions into the recursive form, we have \begin{equation}
    \widehat{\mclk}_2(x,y) = \frac{\mclk(x,s_1)\mclk(t_1,y)}{a} + \frac{ \left(a\mclk(x,s_2)-b\mclk(x,s_1)\right) \left(a\mclk(t_2,y)-c\mclk(t_1,y)\right) }{a(ad-bc)}.
\end{equation}
After simplification, we get  $\widehat{\mclk}_2(x,y) = \dfrac{ (d\mclk(x,s_1)-c\mclk(x,s_2))\mclk(t_1,y) + (-b\mclk(x,s_1)+a\mclk(x,s_2))\mclk(t_2,y) } {ad-bc}.$ 

On the other hand, the compact cross-interpolatory form is $\widehat{\mclk}_2(x,y) = \begin{bmatrix} \mclk(x,s_1) & \mclk(x,s_2) \end{bmatrix} \mclk(T,S)^{-1} \begin{bmatrix} \mclk(t_1,y)\\ \mclk(t_2,y) \end{bmatrix}.$ Using the inverse of $\mclk(T,S)$, we obtain \begin{equation}
    \widehat{\mclk}_2(x,y) = \frac{ (d\mclk(x,s_1)-c\mclk(x,s_2))\mclk(t_1,y) + (-b\mclk(x,s_1)+a\mclk(x,s_2))\mclk(t_2,y) }{ad-bc},
\end{equation}
which coincides exactly with the expression obtained from the recursive update. Hence, the recursive greedy construction and the compact cross-interpolatory representation hold for $k=2$.

\section{Computational Complexity of \autoref{alg: optimal-rank-r-residual-update1}}

In this section, we provide a qualitative analysis of the computational complexity of \autoref{alg: optimal-rank-r-residual-update1}. The dominant cost arises from the optimization step used to determine the pivot pairs $(t_k,s_k)$ at each iteration, while the subsequent kernel updates are comparatively inexpensive.

\begin{itemize}
    \item At each iteration $k$, the pivot pair $(t_k,s_k)$ is obtained by solving the nonlinear optimization problem, i.e., minimization of the residual functional $E_k(t,s)$ given in \eqref{equ: residual functional at kth srep}. The integral in \eqref{equ: residual functional at kth srep} is approximated using quadrature on discretization sets $X_\mathrm{quad} \subset \mclx$ and $Y_\mathrm{quad} \subset \mcly$ of size $N_\mathrm{quad}$, leading to an evaluation cost of $\mclo{N_{\mathrm{quad}}^2}$ per function evaluation.

    \item The cost for evaluating the residual kernel at the quadrature nodes at $k$-th iteration is of $\mclo{2^{2k}\, N_\mathrm{quad}^2}$.

    \item If $N_{\mathrm{init}}$ number of initial guesses are used and each optimization requires $N_{\mathrm{iter}}$ iterations, the total cost of the minimization step per rank is of $\mclo{2^{2k}\, N_{\mathrm{init}} \, N_{\mathrm{iter}} \, N_{\mathrm{quad}}^2}$.


\end{itemize}

Over $r$ iterations, the overall computational complexity is dominated by the repeated optimization and can be estimated as $\mclo{r\,2^{2k}\, N_{\mathrm{init}} \, N_{\mathrm{iter}} \, N_{\mathrm{quad}}^2}$.

The above estimate highlights that the primary computational cost lies in the nonlinear minimization required for node selection. On the other hand, the surrogate residual energy evaluation is relatively inexpensive. In practice, this cost can be minimized by parallelizing across initial guesses and by efficient quadrature implementations.


{\footnotesize
\bibliographystyle{siam}
\bibliography{ref}}
\end{document}